\documentclass[11pt]{article}
\usepackage[T1]{fontenc}
\usepackage{geometry}			
\usepackage[parfill]{parskip}	
\usepackage{csquotes}			
\usepackage{mathtools}		
	\allowdisplaybreaks		
\usepackage{xcolor}
\usepackage{graphicx}			
\usepackage{amsthm}			
\usepackage{amssymb}

\newtheorem{thm}{Theorem}
\newtheorem{prob}[thm]{Problem}

\DeclareMathOperator{\Ber}{Ber}
\DeclareMathOperator{\diag}{diag}
\DeclareMathOperator{\Med}{Med}
\newcommand{\defining}[1]{\emph{#1}}		
\newcommand{\Iff}{if\textcompwordmark f{}}	
\DeclareMathOperator{\fchar}{char}
\DeclareMathOperator{\rank}{rank}
\DeclareMathOperator{\Tor}{Tor}
\newcommand{\defn}{\coloneqq}

\newcommand{\f}{\mathcal{F}}
\newcommand{\M}{\mathcal{M}}
\newcommand{\bC}{\mathbb{C}}
\newcommand{\bP}{\mathbb{P}}
\newcommand{\bR}{\mathbb{R}}
\newcommand{\bF}{\mathbb{F}}
\newcommand{\bfa}{\mathbf{a}}
\newcommand{\bfx}{\mathbf{x}}
\newcommand{\transpose}[1]{#1^{T}}
\newcommand{\calD}{\mathcal{D}}
\newcommand{\calL}{\mathcal{L}}

\usepackage[backref=page]{hyperref}
	\hypersetup{
		colorlinks=true,
		linkcolor=red,
		urlcolor=blue,
		citecolor=magenta
	}
\usepackage[abbrev,lite]{amsrefs}
\usepackage{braha-amsrefs}

\title{An ensemble of high rank matrices arising from tournaments\footnote{\copyright 2022. This manuscript version is made available under the CC BY-NC-ND 4.0 license \url{https://creativecommons.org/licenses/by-nc-nd/4.0/}. The published journal article is available at \emph{Linear Algebra Appl.}\ \textbf{658} (2023), 310--318, \href{https://doi.org/10.1016/j.laa.2022.11.004}{doi:10.1016/j.laa.2022.11.004}.}}
\author{
	Niranjan Balachandran\thanks{Department of Mathematics, Indian Institute of Technology Bombay, Mumbai, India. Email: \href{mailto:niranj@math.iitb.ac.in}{niranj@math.iitb.ac.in}. ORCiD: \url{https://orcid.org/0000-0002-4142-3857}}
	\and
	Srimanta Bhattacharya\thanks{Department of Computer Science and Enginnering, Indian Institute of Technology Palakkad, Palakkad, India. Email: \href{mailto:srimanta@iitpkd.ac.in}{srimanta@iitpkd.ac.in}. ORCiD: \url{https://orcid.org/0000-0002-2612-0217}. This work was done while the author was at the Indian Statistical Institute, Kolkata, India.}
	\and
	Brahadeesh Sankarnarayanan\thanks{Corresponding author. Department of Mathematics, Indian Institute of Technology Bombay, Mumbai, India. Email: \href{mailto:bs@math.iitb.ac.in}{bs@math.iitb.ac.in}. ORCiD: \url{https://orcid.org/0000-0001-9191-1253}. The author's research is supported by the National Board for Higher Mathematics (NBHM), Department of Atomic Energy (DAE), Govt.\ of India.}
}
\date{November 11, 2022}
\begin{document}
\maketitle

\begin{abstract}
Suppose $\mathbb{F}$ is a field  and let $\mathbf{a} \coloneqq (a_1,a_2,\dotsc)$ be a sequence of non-zero elements in $\mathbb{F}$.
For $\mathbf{a}_n \coloneqq (a_1, \dotsc, a_n)$, we consider the family $\mathcal{M}_n(\mathbf{a})$ of $n\times n$ symmetric matrices $M$ over $\mathbb{F}$ with all diagonal entries zero and the $(i, j)$th element of $M$ either $a_i$ or $a_j$ for $i < j$.
In this short paper, we show that all matrices in a certain subclass of $\mathcal{M}_n(\mathbf{a})$---which can be naturally associated with transitive tournaments---have rank at least $\lfloor 2n/3 \rfloor - 1$.
We also show that if $\operatorname{char}(\mathbb{F})\neq 2$ and $M$ is a matrix chosen uniformly at random from $\mathcal{M}_n(\mathbf{a})$, then with high probability $\operatorname{rank}(M) \geq \bigl(\frac{1}{2} - o(1)\bigr)n$.
\end{abstract}

2020 Mathematics Subject Classification: 15B52, 05D99, 05C20, 15A03.\\
Key words and phrases: rank, symmetric matrix, tournament, Talagrand's inequality, bisection closed family

\section{Introduction}
By $[n]$ we shall mean the set $\{1,\dotsc,n\}$.
Suppose $\bF$ is a field and suppose $\bfa \defn (a_1,a_2,\dotsc)$ is a sequence of non-zero elements in $\bF$.
Write $\bfa_n \defn (a_1,\dotsc,a_n)$.
This paper concerns itself with the following 
\begin{prob}\label{probl} 
	Let $\M_n(\bfa)$ consist of the family of all symmetric $n \times n$ matrices over $\bF$ with all diagonal entries being zero and such that for $1\le i < j\le n$ the $(i,j)$th entry is either $a_i$ or $a_j$.
	Determine $\min_{M\in\M_n(\bfa)} \rank(M)$.
	Here, $\rank(M)$ denotes the rank of $M$ over the field $\bF$.
\end{prob}

This problem first appeared in \cite{BalMatMis19}, though it was only stated for matrices over the reals.
There it is asked whether one can find an absolute constant $c>0$ such that $\rank(M)\ge cn$ for all $M\in\M_n(\bfa)$. 

The problem of determining the rank of specific matrices has been of immense interest in extremal combinatorics with applications in theoretical computer science as well---see \cites{Alo15,BarDviWig11,Buk17,CojErgGao20,FalHog07,HuJohDav16,KisSch17,RecFazPar10}.
The question in \cite{BalMatMis19} is motivated by a problem in extremal combinatorics concerning what are called \defining{bisection closed families}: a family $\f$ of subsets of $[n]$ is called a bisection closed family if, for any distinct $A,B\in\f$, either $\frac{\lvert A\cap B\rvert}{\lvert A\rvert}=\frac{1}{2}$ or $\frac{\lvert A\cap B\rvert}{\lvert B\rvert}=\frac{1}{2}$, and one seeks to find the maximum size of a bisection closed family over $[n]$.
One of the results that appears in \cite{BalMatMis19} shows that any bisection closed family has size $O(n\log_{2} n)$, while there are bisection closed families of size $\Omega(n)$.
If the answer to the question in \cite{BalMatMis19} is affirmative, then in fact it is not hard to see that any bisection closed family has size $O(n)$, and we include that simple argument here for the sake of completeness.
Suppose $\f$ is a bisection closed family over $[n]$ of size $m$ and let $X_{m\times n}$ be the matrix whose rows are indexed by the members of $\f$ and the columns by the elements of $[n]$ defined as follows: for $A\in\f$ and $x\in [n]$, set $X(A,x) \defn 1$ if $x\in A$ and $X(A,x) \defn -1$ otherwise.
For two sets $A,B\in\f$, let $\Tor(A,B)\defn A$ if $\lvert A\cap B\rvert=\frac12\lvert B\rvert$ and $\Tor(A,B)\defn B$ otherwise.
Then the matrix $X\transpose{X}$ whose rows and columns are indexed by the members of $\f$ satisfies 
\begin{align*}
	X\transpose{X}(A,A) &{}= n,\\ 
	X\transpose{X}(A,B) &{}= n-2(\lvert A\rvert +\lvert B\rvert)+4\lvert A\cap B\rvert \\
			  			&{}= n-2\lvert \Tor(A,B)\rvert
\end{align*}
In particular, if $\f=\{A_1,\dotsc,A_m\}$ and we let $J$ denote the $m\times m$ matrix consisting entirely of ones, then $\frac12(nJ-X\transpose{X})\in\M_n(\bfa)$ for the sequence $\bfa_n=(\lvert A_1\rvert,\dotsc,\lvert A_m\rvert)$.
Hence, if the conjecture holds, then $\rank(X\transpose{X})\ge cm$. But, since $\rank(X\transpose{X})\le\rank(X) \leq n$, it follows that $m \leq (n + 1)/c$, which establishes an asymptotically tight bound on the size of the bisection closed family.

In this short paper, we make some steps towards settling this problem in the affirmative.
In order to describe our results, we note that to each $M\in\M_n(\bfa)$ there corresponds a tournament on the vertex set $[n]$ in the following natural manner: for $i<j$ we direct the edge $ij$ as $i\to j$ if $M(i,j)=a_i$, and the edge is directed in the reverse direction if $M(i,j)=a_j$.
Conversely, for a tournament $T$ on $[n]$, we can associate the matrix $M_T(\bfa)\in\M_n(\bfa)$ in exactly the same way, namely, for $i<j$, set $M_T(i,j)=a_i$ \Iff\ $i\to j$.
Note that this correspondence is not necessarily one-to-one, since the $a_i$ need not be distinct. 

Our first result gives a lower bound in the case where the underlying tournament is transitive, i.e., when there is a total order \(\prec\) on \([n]\) such that $i\to j$ whenever \(i \prec j\).

\begin{thm}\label{trans_tour}
	If $T$ is transitive, then $\rank(M_T(\bfa))\ge \bigl\lfloor \frac{2n}{3} \bigr\rfloor - 1$.
\end{thm}

In \cite{BalMatMis19}, there are constructions of bisection closed families $\f$ of size $\frac{3}{2}n-2$ which admit a uniform subfamily of size $n - 1$. So, by the remarks above, the corresponding matrix $M$ has a principal submatrix of rank at least $2n/3 + 1$, so the constant $2/3$ in the theorem is best possible over \emph{all} tournaments. It is instructive to compare this with a result of de Caen's \cite{Cae91} on the rank of tournament matrices (where the entries are only $0$ and $1$). In \cite{Cae91}, among other things, it is shown that the rank of any $n \times n$ tournament matrix is at least $\frac{n-1}{2}$ over any field and at least $n-1$ over the reals. Our result is in a similar spirit, since it answers a question in a more general setup, but is also in contrast with de Caen's since we consider symmetric matrices.

Our second result shows that almost all the matrices in $\M_{n}(\bfa)$ have high rank.
More precisely, we show that for a \defining{random tournament}---a tournament with the edges being directed in either direction with probability $1/2$ each and independently---then \defining{with high probability (whp)} the rank is at least $(1/2-o(1))n$. Here, the phrase \enquote{with high probability} means that the probability that the said event occurs asymptotically tends to $1$ as $n\to\infty$.

\begin{thm}\label{rand_tour}
Suppose $\fchar(\bF)\ne 2$ and $\bfa$ is a sequence of non-zero elements of $\bF$. If $T$ is a uniformly random tournament, then \emph{whp} $\rank(M_T(\bfa))\ge \frac{n}{2}-21\sqrt{n\log n}$.
\end{thm}

To give some perspective on this result vis-\`a-vis the existing literature on similar problems, the behavior of random symmetric matrices is an immensely active area of research and there are several papers that consider various random models (see \cites{FerJai19,TaoVu17} and the references therein) and Theorem~\ref{rand_tour} may be regarded as another addition to that list, though there is a fundamental difference between our result and all the others.
For one, as we have pointed out earlier, the matrix that arises from a bisection closed family has rank at most $2n/3+O(1)$, so in that sense our result is somewhat qualitatively different from those that appear in several of those papers.
It must be pointed out that the main result in \cite{TaoVu17} considers random symmetric matrices $M_n=((\xi_{ij}))$ where $\xi_{ij}$  are all jointly independent (for $i<j$) and also independent of $\xi_{ii}$ (which are also independent) with the additional property that for all $i<j$ and all real $x$, $\bP(\xi_{ij}=x)\le 1-\mu$ for some fixed constant $\mu$, and their result shows that \emph{whp} the spectrum is simple.
This does establish (in a strong form) Theorem~\ref{rand_tour} over the reals, {\it in the special case where $a_i$ are all pairwise distinct}.
But otherwise the best bound this suggests is of the order $\Omega(\sqrt{n})$.
Secondly, our result holds over all fields $\bF$ with $\fchar(\bF)\ne 2$ whereas most other results usually work specifically with $\bR$ or $\bC$ (though they have stronger results).
To also contrast the results of Theorems~\ref{trans_tour} and \ref{rand_tour}, note that Theorem~\ref{trans_tour} holds over all fields whereas for Theorem~\ref{rand_tour} we need $\fchar(\bF)\ne 2$.

We prove Theorems~\ref{trans_tour} and \ref{rand_tour} in the next section.
The final section includes some concluding remarks and poses some further questions. 

\section{Proofs of Theorems~\ref{trans_tour} and \ref{rand_tour}}

\begin{proof}[Proof of Theorem~\ref{trans_tour}]
It suffices to prove that whenever \(3\) divides \(n\) we have \(\rank(M_{T}(\bfa)) \geq 2n/3\), since we may then interpolate to those \(n\) such that \(3\) does not divide \(n\) to show that \(\rank(M_{T}(\bfa)) \geq \lfloor 2n/3 \rfloor - 1\) for all \(n \geq 3\).
We shall prove this by induction on $n \geq 3$ such that \(3\) divides \(n\).
Without loss of generality we may assume that the ordering of the elements coincides with the natural ranking order on $[n]$, i.e., $i\prec j$ \Iff\ $i>j$.
Also, we denote the matrix corresponding to the transitive tournament on $[n]$ by $D_n(a_{1},\dotsc,a_{n})$.
When the $a_{i}$ are clear from the context, we will simply call this matrix $D_{n}$.
Then, for the base case, i.e., for $n=3$, it is easy to see that 
\[
	D_3(a_{1},a_{2},a_{3})=
	\begin{pmatrix}
		0     & a_{2} & a_{3} \\
		a_{2} & 0     & a_{3} \\
		a_{3} & a_{3} & 0	
	\end{pmatrix}
\]
has rank at least \(2 = \frac{2 \cdot 3}{3}\) for any non-zero values of $a_2$ and $a_3$.
Next, we assume the assertion to be true for $n$ and prove it for $n+3$.
Then, up to relabeling of the indices, we can write $D_{n+3}(a_{1},\dotsc,a_{n+3})$ as
\begin{equation}\label{E:block}
	D_{n+3}=
	\begin{pmatrix}
		D_{3} 		  & B \\
		\transpose{B} & D_{n}
	\end{pmatrix},
\end{equation}
where $D_{3} = D_{3}(a_{1},a_{2},a_{3})$, $D_{n} = D_{n}(a_{4},\dotsc,a_{n+3})$,
and $B$ is the matrix 
\[
	\begin{pmatrix}
		a_4 & a_5 & \cdots & a_{n+3}\\
		a_4 & a_5 & \cdots & a_{n+3}\\
		a_4 & a_5 & \cdots & a_{n+3}
\end{pmatrix}
\]
and by the induction hypothesis $D_{n}$ has rank $2n/3$.
Now, assume that $D_{n+3}$ has rank $\leq \frac{2n}{3}+1$.
Let \(\calL\) be a basis of $D_{n+3}$.
So, by our assumption $\lvert \calL \rvert \leq \frac{2n}{3}+1$.
Next, we will arrive at a contradiction by considering the following cases.

For ease of presentation, let us write $D_{n+3} = \begin{pmatrix} \bfa_1 & \cdots & \bfa_n \end{pmatrix}$, $D_{3} = \begin{pmatrix} \bfa'_{1} & \bfa'_{2} & \bfa'_{3} \end{pmatrix}$, and $\transpose{B} = \begin{pmatrix}\bfa''_{1} & \bfa''_{2} & \bfa''_{3} \end{pmatrix}$.
Let $\calD = \{\bfa_1, \bfa_2, \bfa_3\}$.

\textsf{Case 1:} $\lvert \calL \cap \calD\rvert =0$.
In this case, $\calL \subseteq \{\bfa_4, \dotsc, \bfa_n\}$.
But, this is not possible since the column space of $B$, which is of rank one, cannot include the column space of $D_{3}$, which is of rank three.
In fact, no column of $D_{3}$ belongs to the column space of matrix $B$.
More concretely, for $\bfa_i \in \calD, 1 \leq i \leq 3$, let $\bfa_i = \sum\gamma_j \bfa_j$, where $\gamma_j \in \bF$ and $\bfa_j \in \calL$ with $j \in \{4, \dotsc, n+3\}$.
Then, it follows that $\bfa'_i = \sum \gamma_j \bfa'_j$, where $\bfa'_j = \transpose{\begin{pmatrix} a_j & a_j & a_j \end{pmatrix}}$.
But, this is a contradiction since $0= \sum_j \gamma_j a_{j} = a_3 \neq 0$.

\textsf{Case 2:} $\lvert \calL \cap \calD\rvert =1$.
This is similar to \textsf{Case 1}.
More precisely and without loss of generality, let $\calL$ contain $\bfa_3$.
Now, let $\bfa_2 = \gamma_3 \bfa_3 + \sum_j \gamma_j \bfa_j$, where $\gamma_j \in \bF$ and $\bfa_j \in \calL$ with $j \in \{ 4, \ldots, n+3\}$.
This implies $\bfa'_{2} = \gamma_3 \bfa'_{3} + \sum_j \gamma_j \bfa'_{j}$.
Then, following the considerations of \textsf{Claim 1} we have $\gamma_3 \neq 0$.
But then we have that $0 = \gamma_3 a_3 + \sum_j \gamma_j a_j = a_2 \neq 0$ for $j \in \{4, \dotsc, n+3\}$, which is a contradiction.
Similar arguments hold for the cases when $\bfa_1 \in \calL$ or $\bfa_2 \in \calL$.

\textsf{Case 3:} $\lvert \calL \cap \calD\rvert \in \{2, 3\}$ for any choice of \(\calL\).
This implies that $\rank\bigl(\transpose{\begin{pmatrix}\transpose{B} & D_{n} \end{pmatrix}}\bigr) \leq \frac{2n}{3} - 1$, which is a contradiction since $\rank(D_{n})$ is already at least $2n/3$ by our assumption.
\end{proof}

Before we get to the proof of Theorem~\ref{rand_tour}, we state a few results that we shall use.
We state the versions as they appear in \cite{MolRee02}.
\begin{thm}[Chernoff bound]
	If $X$ is distributed as the binomial random variable $B(n,p)$, then for any $0\le t\le np$
	\[
		\bP(\lvert X-np \rvert>t) < 2\exp\left(-\frac{t^2}{3np}\right).
	\]
\end{thm}
The other main technical tool is Talagrand's inequality.
By a \emph{trial} we shall simply mean a Bernoulli event. 
\begin{thm}[Talagrand's inequality]\label{T:Talagrand}
	Suppose $X$ is a non-negative random variable, not identically zero, which is determined by $n$ independent trials $T_1,\ldots, T_n$, and satisfying the following for some $c,r>0$:
	\begin{enumerate}
		\item (\textbf{$c$-Lipschitz}) Changing the outcome of any one trial $T_i$ changes $X$ by at most $c$,
		\item (\textbf{$r$-certifiability}) For any $s\ge 0$, if $X\ge s$ then there is a set of at most $rs$ trials whose outcomes certify that $X\ge s$, i.e., there is a set $I\subset [n]$ of size at most $rs$ and a set of outcomes of the trials $T_i$ for $i\in I$ such that  fixing the outcomes of $T_i$ for $i\in I$ ensures $X\ge s$  \emph{irrespective of the outcomes of $T_j$} for $j\notin I$.
	\end{enumerate}
	If $\Med(X)$ denotes the median of $X$ and $0\le t\le \Med(X)$, we have
	\[
		\bP(\lvert X-\Med(X)\rvert>t)\le 4\exp\left(-\frac{t^2}{8c^2r\Med(X)}\right).
	\]
\end{thm}

\begin{proof}[Proof of Theorem \ref{rand_tour}]
We begin with a couple of simple observations. 
\begin{enumerate}
	\item Fix a pair $(i,j)$ with $i<j$, and let $T$ be a tournament on $[n]$.
	If $T'$ is the  tournament obtained from $T$ by changing the orientation of only the edge $ij$, then $M_{T'}(\bfa)=M_T(\bfa)+D$ for a matrix $D$ comprising of zeros everywhere except at the $(i,j)$ and $(j,i)$ positions.
	Consequently, $\lvert\rank(M_T(\bfa))-\rank(M_{T'}(\bfa))\rvert\le 2$.
	Also, observe that if $\bfa'$ is the sequence with $a_i$ replaced by some $z\in\bF$, then for any tournament $T$, $M_T(\bfa)$ differs from $M_T(\bfa')$ only in the entries of the \(i\)th row and column, so again in particular, $\lvert \rank(M_T(\bfa))-\rank(M_T(\bfa'))\rvert\le 2$.
	
	\item For a tournament $T$, let $T_R$ denote the \defining{reverse} tournament, i.e., if $i\to j$ in $T$ then $j\to i$ in $T_R$.
	Then $M_T(\bfa)+M_{T_R}(\bfa)=M$ where $M(i,i)=0$ and $M(i,j)=a_i+a_j$.
	In particular, $M=DJ+JD-2D$ where $D$ is the diagonal matrix $\diag(a_1,\ldots,a_n)$ and $J$ represents, as before, the all-ones matrix.
	In particular, since $a_i\ne 0$ and \(\fchar(\bF) \neq 2\), it follows that $\rank(M)\ge n-2$.
	Consequently, at least one of $\rank(M_T)$ and $\rank(M_{T_R})$ is at least $n/2 - 1$.
\end{enumerate}

First, suppose that $\fchar(\bF)$ does not divide $n-1$.
A uniformly random tournament $T$ is completely determined by the trials $T_e$ for the pairs $e=(i,j)$ with $i<j$ with $T_e$ distributed as $\Ber(1/2)$.
Fix some $z\ne 0$ in $\bF$ and consider the more general ensemble $\M_n(\bfx)$, where $\bfx=(x_1,\ldots,x_n)$ is the random sequence with $x_i=z$ with probability $1/\sqrt{n}$ and $x_i=a_i$ with probability $1-\frac{1}{\sqrt{n}}$.
(We can view \(\bfx\) as arising from \(n\) flips of a biased coin where heads occurs with probability \(1/\sqrt{n}\), and \(x_{i} = z\) if a head occurs on the \(i\)th toss
and \(x_{i} = a_{i}\) if a tail occurs on the \(i\)th toss.)
Let $X=\rank(M_T(\bfx))$.
By the first observation, it follows that $X$ is $2$-Lipschitz.
Also, for any $s\ge 1$, if we fix $x_1=\dotsb =x_s= x_{s+1} =z$, then, irrespective of $T$, either the principal $s \times s$ matrix or the principal \((s+1) \times (s+1)\) submatrix of $M_T(\bfx)$ has non-zero determinant in $\bF$, so this establishes that $X$ is $2$-certifiable.
Note that if $T$ is a transitive tournament, then so is $T_R$, so it follows by Theorem~\ref{trans_tour} and the second observation above that $\Med(X)\ge \frac{n-2}{2}$.
Hence, by Talagrand's inequality,
\[
	\bP\left(X > \frac{n}{2} - 16\sqrt{n\log n}\right)<\frac{4}{n^4}.
\]
Now, let $Y$ denote the number of elements of $\bfx$ that are equal to $z$ (that is, the number of occurrences of heads in the sequence of coin flips).
By the Chernoff bound, it follows that $\bP(Y>2\sqrt{n})<2e^{-\sqrt{n}/3}$, so with probability at least $1-\frac{4}{n^4}-2e^{-\sqrt{n}/3}$ both the events, $X>n/2 - 16\sqrt{n\log n}$ and $Y<2\sqrt{n}$ hold simultaneously.
Hence, if we only consider the submatrix $M'$ of $M_T(\bfx)$ indexed by those rows and columns of the sequence $\bfx$ that do not comprise of any $z$, then $\rank(M')\ge n/2 - 16\sqrt{n\log n}-4\sqrt{n}$.
Now finally, if $\fchar(\bF)$ divides $n-1$, then we restrict ourselves to the principal submatrix of $M_T$ of order $n-1$, and the same argument as above gives us $\rank(M_T(\bfa))\ge n/2 - 21\sqrt{n\log n}$ \emph{whp} with room to spare.
\end{proof}

\section{Concluding remarks}\label{S:conclusion}

\begin{itemize}
	\item In the statement of Theorem~\ref{rand_tour}, the randomness is over the orientations of the tournament edges, and \emph{not} over the elements $a_i$ since they come from the given sequence $\bfa$.
	However, the proof uses a randomization of the sequence in order to be able to use Talagrand's inequality, and it does not seem straightforward to stay within the confines of the given sequence to be able to prove the same statement. 

	\item Our bound of $n/2$ is constricted by our estimate of $\Med(\rank(M_T(\bfa))$.
	If one can get a better bound, then the same proof gives a better rank bound as well.
	Interestingly, this is an instance where using Talagrand's inequality for concentration around the median gives a decidedly better bound than the analogous version for concentration around the mean.

	\item Our error probability of $O(n^{-3})$ is easily improved to $O\bigl(e^{-n^{1/3}}\bigr)$ if we take $t=n^{2/3}$ for instance in the proof of Theorem~\ref{rand_tour}.
	
	\item As remarked after the statement of Theorem~\ref{trans_tour}, while the lower bound in Theorem~\ref{trans_tour} does achieve the bound that is best possible for all tournaments, we believe that the bound must be substantially better when restricted to transitive tournaments.
	In fact, we believe that for transitive tournaments, $\rank(M_T(\bfa))\ge n-o(n)$ must hold as well though we are unable to prove this even over the reals.\footnote{This has recently been settled in the affirmative in a strong form; see~\cite{BalBhaSan22}.}

	\item As indicated in the remarks in the introduction, it must be possible to improve upon the results obtained here when we restrict ourselves to the fields $\bR$ or $\bC$, or if the elements $a_i$ themselves satisfy other constraints.
	For instance, when the sequence $\bfa$ is the constant sequence $(a,a,\dotsc)$, $\rank(M) \geq n - 1$ for all $M \in \M_{n}(\bfa)$ and over any field. This also shows that if $\bF$ is any finite field, then $\rank(M) \geq \frac{n}{\lvert \bF \rvert - 1} - 1$ for all $M \in \M_{n}(\bfa)$ for \emph{any} sequence $\bfa$ in $\bF$. On the other hand, it is not even clear whether the results of this paper extend to the case when $\fchar(\bF) = 2$, or when infinitely many of the $a_{i}$ are distinct.
	
	\item A more general setup is the following: Suppose $f\colon \bF^2\to\bF$ and let $\bfa=(a_1,a_2,\dotsc)$ as before and consider the central problem of this paper over the more general family $\M^{(f)}_{n}(\bfa)$ which is defined as follows.
	For any tournament $T$ on the vertex set $[n]$ the matrix $M_T^{(f)}(\bfa)\in\M^{(f)}_{n}(\bfa)$ consists of zeros on the diagonal, and for $i<j$ the $(i,j)$ entry of $M_T^{(f)}(\bfa)$ equals $f(a_i,a_j)$ if $i\to j$ in $T$, and equals $f(a_j,a_i)$ otherwise.
	Our proof of Theorem~\ref{rand_tour} is easily modified to show that for a random tournament $T$, $\rank\bigl(M_{T}^{(f)}(\bfa)\bigr) \geq \bigl(\frac12 - o(1)\bigr)n$ for any function $f(x,y)=\alpha x+\beta y$ with the property $\alpha+\beta\ne 0$. The only difference in the proof is that we use Talagrand's inequality with concentration about the mean instead of the median.
	If $f(x,y)=\alpha x+\beta y$ with $\alpha+\beta\ne 0$, then the same argument as in the proof of Theorem~\ref{rand_tour} also shows that $\rank\bigl(M_T^{(f)}(\bfa)\bigr)+\rank\bigl(M_{T_R}^{(f)}(\bfa)\bigr)\ge n-2$, so again for a random tournament $T$ the expected rank of $M_T^{(f)}(\bfa)$ is at least $n/2 -1$. We omit the details.\footnote{Note that resolving the problem for the family $\M_{n}^{(f)}(\bfa)$ also gives a linear upper bound on the size of any $\theta$-intersecting family (cf.~\cite{BalMatMis19}), in a similar manner.}
	
	This more general ensemble of matrices may pose yet more interesting difficulties even for relatively simple functions. For instance, even for $f(x,y)=\frac{x}{y}$, the problem is already quite non-trivial.
	
	\item For a given sequence $\bfa$ in $\bF$ and a matrix $M = M_T(\bfa_{n}) \in \M_{n}(\bfa)$ arising from a self-dual tournament $T$, the matrix $M_{T_{R}}(\bfa_{n})$ arising from the reverse tournament $T_{R}$ can also be viewed as $M_{T}(\sigma \bfa_{n})$ for a permutation $\sigma$ of $\bfa_{n} = (a_{1},\dotsc,a_{n})$.
	We have shown that at least one of $M_{T}$ and $M_{T_{R}}$ has rank at least $n/2 - 1$ for any tournament $T$ on $[n]$.
	An interesting question is whether, for a fixed tournament $T$, the matrices $M_{T}(\sigma \bfa_{n})$ have the same rank for all permutations $\sigma$ of $\bfa_{n}$.
	A positive answer to this question will tell us, in particular, that matrices arising from self-dual tournaments (such as Paley tournaments) all have high rank.
\end{itemize}

\appendix
\newpage
\section*{ADDENDUM\let\thefootnote\relax\footnote{Date: July 13, 2023}}

We fix a small error in the proof of Theorem~\ref{rand_tour}.
The statement of the theorem remains unchanged.

\subsection*{The proof of Theorem~\ref{rand_tour}}
Our application of Talagrand's inequality (Theorem~\ref{T:Talagrand}) in the proof of Theorem~\ref{rand_tour} is incorrect: it is true that if we fix $x_1 = \dotsb = x_s = x_{s+1} = z$ gives $X \geq s$, but this is not sufficient to establish that $X$ is $2$-certifiable.
Instead, we need a version of McDiarmid's inequality for concentration bounds on product measure spaces, and we use the one stated in \cite[Lemma 1.2]{McDiarmid}:
\begin{thm}[Independent Bounded Differences Inequality]
    Let $X_{1},\dotsc,X_{n}$ be independent random variables, with $X_{k}$ taking values in a set $\Omega_{k}$ for each $k$.
    Suppose that the measurable function $f \colon \prod_{k} \Omega_{k} \to \mathbb{R}$ satisfies, for each $k$,
    \[
    	\lvert f(\mathbf{x}) - f(\mathbf{x}') \rvert \leq c_{k}
    \]
    whenever the vectors $\mathbf{x}$ and $\mathbf{x}'$ differ only in the $k$th coordinate.
    Let $Y$ be the random variable $f(X_{1},\dotsc,X_{n})$.
    Then, for any $t > 0$,
    \[
    	\mathbb{P}(\lvert Y - \mathbf{E}(Y) \rvert > t)\, \leq\, 2 \exp\bigl(-2t^{2} / {\textstyle{\sum_{k}}} c_{k}^{2}\bigr).
    \]
\end{thm}

\begin{proof}[Proof of Theorem~\ref{rand_tour}]
	The notation $[m,n]$ denotes the set of integers $i$ such that $m \leq i \leq n$.
	Let $\Omega_{k} = \{ 0, 1 \}^{[k+1,n]}$ for $k = 1,\dotsc,n-1$.
	View each vector $\mathbf{x}_{k} = (x_{k}^{(k+1)},\dotsc,x_{k}^{(n)}) \in \Omega_{k}$ as a win-loss record for player $k$ against the players $k+1,\dotsc,n$ in that order.
	Thus, any $(n-1)$-tuple $(\mathbf{x}_{1},\dotsc,\mathbf{x}_{n-1})$, where $\mathbf{x}_{k} \in \Omega_{k}$ for each $k$, determines a unique tournament $T$ on $[n]$, and each tournament on $[n]$ arises from some point in $\prod_{k} \Omega_{k}$.
	
	Define $f \colon \prod_{k} \Omega_{k} \to \mathbb{R}$ by $f(\mathbf{x}_{1},\dotsc,\mathbf{x}_{n-1}) = \operatorname{rank}(M_{T}(\mathbf{a}))$, where $T$ is the tournament uniquely determined by $\mathbf{x}_{1},\dotsc,\mathbf{x}_{n-1}$.
	As observed previously, changing the orientation of any edge $ij$ to get a tournament $T'$ changes the rank by at most $2$, since the matrix $M_{T'}(\mathbf{a})$ is obtained from $M_T(\mathbf{a})$ by adding a matrix of rank at most $2$.
	In fact, for any fixed $i$, flipping the orientations of any subcollection of the edges $ij$, where $j \geq i$, changes the rank by at most $2$, since this again corresponds to adding a matrix of rank at most $2$ to $M_{T}(\mathbf{a})$.
	Hence, $\lvert f(\mathbf{x}) - f(\mathbf{x}') \rvert \leq 2$ whenever $\mathbf{x}$ and $\mathbf{x}'$ differ only in the $k$th coordinate.
	
	Now, let $T$ be a uniformly random tournament on $[n]$, i.e., one for which the orientation of each edge $ij$ is chosen by a fair coin toss.
	Then, $T$ gives rise to random variables $X_{k}$ taking values in $\Omega_{k}$ for each $1 \leq k \leq n - 1$, and $X_{1},\dotsc,X_{n-1}$ are independent.
    Define $Y \coloneq f(X_{1},\dotsc,X_{n-1})$.
    As observed previously, if $T_R$ denotes the reverse tournament of $T$, then $\operatorname{rank}(M_T(\mathbf{a}) + M_{T_R}(\mathbf{a})) \geq n - 2$, so $\mathbf{E}(Y) \geq \frac{n}{2}-1$.

	Now, we apply McDiarmid's inequality to get:
    \[
        \mathbb{P}\left(Y < \frac{n}{2} - 1 - 4\sqrt{n \log n}\right) \leq 2 e^{-\frac{32 n \log n}{4(n-1)}} < \frac{2}{n^8},
    \]
    which proves the result.
\end{proof}

The above proof is simpler than the original attempt in that we do not randomize the sequence $\mathbf{a}$ in addition to picking the tournament $T$ on $[n]$ at random.
Also, the above proof uses concentration around the mean, and not around the median as in Talagrand's inequality.
Consequently, our first two comments in Section~\ref{S:conclusion} are no longer relevant.

\subsection*{A question on general ensembles of the form $\mathcal{M}_{n}^{(f)}(\mathbf{a})$}

In Section~\ref{S:conclusion}, we defined a more general ensemble $\mathcal{M}_{n}^{(f)}(\mathbf{a})$ as consisting of those symmetric matrices $M_{T}^{(f)}(\mathbf{a})$ with zero diagonal such that the \((i,j)\)th entry is either $f(a_{i},a_{j})$ or $f(a_{j},a_{i})$ depending on the orientation of the edge $ij$ in the tournament $T$.
We raised the problem of finding a lower bound on the rank of the matrices coming from this general ensemble, and observed that the same methods work for linear functions of the form $f(x,y) = \alpha x + \beta y$ for which $\alpha + \beta \neq 0$, and when \(\fchar(\bF) \neq 2\), and that the problem is non-trivial for other (even relatively simple) functions, such as $f(x,y) = x/y$.

As it turns out, this specific example does not illustrate the non-triviality of this question, since the same methods suffice to prove a linear lower bound (\emph{whp}) for functions of finite rank.
More precisely, let \(\fchar(\bF) \neq 2\), and suppose that \(f \colon \bF^{2} \to \bF\) is a function such that \(f(x,x) \neq 0\) for all \(x\).
Suppose that there exist functions \(g_{i}, h_{i} \colon \bF \to \bF\), \(1 \leq i \leq k\), such that \(f(x,y) = \sum_{i=1}^{k} g_{i}(x)h_{i}(y)\).
For \(\mathbf{z} \in \bF^{n}\), define \(G_{i}(\mathbf{z}) \defn \begin{pmatrix} g_{i}(z_{1}) & \cdots & g_{i}(z_{n}) \end{pmatrix}\) and \(H_{i}(\mathbf{z}) \defn \begin{pmatrix} h_{i}(z_{1}) & \cdots & h_{i}(z_{n})\end{pmatrix}\) for all \(1 \leq i \leq k\).
Then, for any tournament \(T\) on \([n]\) and any sequence \(\bfa\) in \(\bF\), we have \(M_{T}^{(f)}(\bfa) + M_{T_{R}}^{(f)}(\bfa) = \sum_{i = 1}^{k} \bigl(G_{i}(\bfa_{n})^{T} H_{i}(\bfa_{n}) + H_{i}(\bfa_{n})^{T} G_{i}(\bfa_{n})\bigr) - 2\diag\bigl(f(a_{1},a_{1}),\dotsc,f(a_{n},a_{n})\bigr)\).
The RHS is the sum of a diagonal matrix---of full rank---and \(2k\) matrices of rank one.
Hence, at least one of \(M_{T}^{(f)}\) or \(M_{T_{R}}^{(f)}\) has rank at least \((n-2k)/2\).
Hence, the proof of Theorem~\ref{rand_tour} above shows that, even in this case, we have \emph{whp} \(\rank\bigl(M_{T}^{(f)}(\bfa)\bigr) \geq \frac{n}{2} - o(n)\) for a uniformly random tournament \(T\) on \([n]\).

It is worth noting that the condition \(f(x,x) \neq 0\) for all \(x\) is crucial. For example, let \(\bF = \mathbb{Q}\), $f(x,y) = (x - y)^{2}$ and $\mathbf{a}_{n} = (1,2,\dotsc,n)$, then $\operatorname{rank}(M) = 3$ for all \(M \in \mathcal{M}_{n}^{(f)}(\bfa)\), $n \geq 3$ \cite[Theorem 2.4]{GroodEtAl2014}.
Thus, we refine our original question and ask: for what functions $f \colon \mathbb{F}^{2} \to \mathbb{F}$ does there exist a constant $c > 0$ such that $\operatorname{rank}(M) \geq cn$ for all $M \in \mathcal{M}_{n}^{(f)}(\mathbf{a})$?

\end{document}